\documentclass[12pt]{article}
\usepackage{amsmath,amssymb}
\usepackage{amsthm}

\begin{document}
\author{Richard D. Burstein}
\title{Commuting-square subfactors and central sequences}
\maketitle

\newcommand{\commsq}[4]{\begin{array}{ccc}#3&\subset&#4\\
\cup& &\cup\\
#1&\subset&#2\end{array}}

\newcommand{\tr}{ {\rm tr} }

\newtheorem{lem}{Lemma}

\newtheorem{thm}{Theorem}

\newtheorem{cor}{Corollary}

\begin{abstract}

Let $M_0 \subset M_1$ be a finite-index infinite-depth
hyperfinite $II_1$ subfactor and $\omega$ a free ultrafilter of
the natural numbers.  We show that if this subfactor
is constructed from a commuting square then the central
sequence inclusion $M_0^{\omega} \cap M_1' \subset
(M_0)_{\omega}$ has infinite Pimnser-Popa index.
We will also
demonstrate this result for certain infinite-depth
hyperfinite
subfactors coming from groups.

\end{abstract}

\section{Introduction}

The central sequence subfactor induced by a finite-depth inclusion
was described in \cite{Oc}, \cite{KaA} and \cite{KaC} (see
also~\cite{EK}).
Let $M_0 \subset M_1 \cong R$ be a finite-depth finite-index hyperfinite subfactor,
with Jones tower $M_0 \subset M_1 \subset M_2 \subset ... \subset M_{\infty}$,
and $\omega$ a free ultrafilter.  In this case the Von Neumann
algebra $M_0^{\omega} \cap M_1'$ is a factor.
It is a finite-index subfactor of $M_1' \cap M_1^{\omega} =
(M_1)_{\omega}$, and the dual of its standard invariant may be
computed from the asymptotic inclusion
$M_0 \vee (M_0' \cap M_{\infty}) \subset M_{\infty}$.

If the original subfactor is infinite depth then the situation is less clear.
$M_0^{\omega} \cap M_1'$ need not be a factor.  We may
still define its index in $(M_1)_{\omega}$ using the generalized
index of Pimsner and Popa \cite{PP}, but it is not known in general whether
this index is finite or infinite.  V. Jones has conjectured
that the index is always infinite in this case, and we will prove this
conjecture for infinite-depth
subfactors constructed in certain ways.

For some infinite-depth subfactors coming from groups, it is possible
to write down $M_0^{\omega} \cap M_1'$ as the fixed points of
the $II_1$ factor $R_{\omega}$ under the action of an infinite
group.  The proof of the conjecture follows immediately in these
cases, as we will show in section 2.
However, it is generally difficult to explicitly describe
$M_0^{\omega} \cap M_1'$; for many examples, it is not clear
whether or not this algebra is larger than $\mathbb{C}$.

In certain other cases, the conjecture can be proved by exhibiting projections
in $(M_1)_{\omega}$ whose conditional expectations onto
$M_0^{\omega} \cap M_1'$ have arbitrarily small norm.  From \cite{PP}, this
gives infinite index for the inclusion, but this approach presents
certain technical difficulties.
If the original subfactor
does not admit a generating tunnel, it may be hard to specify
any elements of $M_0^{\omega} \cap M_1'$ besides the scalars.
While finding many elements of $(M_1)_{\omega}$ is relatively
straightforward, determining their conditional expectations
onto $M_0^{\omega} \cap M_1'$ may be intractable.

Using commuting-square subfactor ameliorates these problems.
Here the tower $M_0 \subset M_1 \subset M_2 \subset ...$
may be approximated by a grid of finite-dimensional
Von Neumann algebras as in \cite{JS}.
As we will show in section 3, the
iterated canonical shift provides projections in $(M_k)_{\omega}$
whose conditional expectation onto
 to $M_0^{\omega} \cap M_1'$ may be bounded in norm.
In section 4 some technical lemmas on index of ultrapower factors,
combined with Sato's result \cite{Sa} on the global index of the horizontal
and vertical
commuting-square subfactors, will give the desired inequalities.

Throughout this paper, we will take $\omega$ to be an arbitrary free
ultrafilter.

\vspace{0.5cm}

\noindent{\it Acknowledgement.}
I am grateful to Professor Vaughan Jones for suggesting this problem to me,
and for outlining the method of solution in the Bisch-Haagerup subfactor case.
The main theorem in this paper was part of my doctoral thesis at
UC Berkeley \cite{RB}.

\section{Subfactors obtained from groups}

\subsection{The diagonal subfactor}

The hyperfinite diagonal subfactor (see \cite{PoF},\cite{BiE})
may be constructed
from a finitely generated group of automorphisms of the
hyperfinite $II_1$ factor $R$.
Let $\alpha_2 ... \alpha_n$
be outer automorphisms of $R$, with $\alpha_1 = id$.  Let
$A=M_n(\mathbb{C})$ have matrix units $\{e_{ij}\}$,
and $M_1 = R \otimes A$.  Each $\alpha_i$ extends to an action
on $M_1$ as $\alpha_i \otimes id$.

Let $\alpha$ be the map from $R$ to $M_1$
given by $\alpha(y) = \sum_{i=1}^n e_{ii} \alpha_i(y)$.  Then
$M_0 = \alpha(R)$,
and  $M_0 \subset M_1$ is the diagonal subfactor.

Let $G$ be the image in ${\rm Out} R = {\rm Aut} R / {\rm Int} R$ of the
group generated by the $\alpha_i$'s, and for $g \in G$
let $\alpha_g$ be an arbitrary representative of $g$ in
${\rm Aut} R$.
All of the $M_0-M_0$
bimodules in the Jones tower
of the subfactor are isomorphic to one of $_{\alpha_g(R)}(R)_{R}$,
for some $g \in G$.  Therefore if the subfactor is infinite depth,
then $G$ is infinite.

Now we compute $M_0^{\omega} \cap M_1'$.  If $x=(x_n) \in (M_1)_{\omega}$
is an
element of $M_0^{\omega} \cap M_1'$, this means that $x_n = \alpha(y_n)$
for some bounded sequence $y=(y_n)$ in $R^{\omega}$.
$x$ asymptotically commutes with the constant
sequence $(1 \otimes e_{i1})$ for each $i$.  Since $(1 \otimes
e_{i1}) x_n =y_n \otimes e_{i1}$,
while $x_n (1 \otimes e_{i1}) = \alpha_i(y_n) \otimes e_{i1}$, it follows that
$||\alpha_i(y_n)-y_n||_2$ goes to zero along the ultrafilter
for all $i$.  From the definition of $\alpha$,
this means that $||\alpha_i(x_n)-x_n||_2$ goes to zero along
the ultrafilter as well.

$\alpha_i$ induces a pointwise
action $(\alpha_i)_{\omega}$ on $(M_1)_{\omega}$ given by
$(\alpha_i)_{\omega}( (a_n) ) = (\alpha_i(a_n))$. 
The above limit then means that
$(\alpha_i)_{\omega}$ fixes $x$, for all $x \in M_0^{\omega} \cap M_1'$.
The same is true of
$(\alpha_g)_{\omega}$ for all $g \in G$, since the
$\alpha_i$'s generate $G$. 

From \cite{Co}, for $\alpha \in {\rm Aut} R$, if $\alpha$ is properly
outer then $\alpha_{\omega}$ is as well.
$M_1 \cong R$, so
$g \rightarrow (\alpha_g)_{\omega}$ defines an outer action
of $G$ on $(M_1)_{\omega}$.  $M_0^{\omega} \cap M_1'$ is contained
in the fixed points of this action.

In \cite{PP}, the authors define
a generalized index on an inclusion of Von Neumann algebras $A \subset B
\subset M$, where $M$ is a $II_1$ factor.  This index
$[B:A]$ is the supremum of $\{\lambda^{-1}\}$, where
$\lambda$ is a scalar such that $E_{A}(x) \ge \lambda x$
for all positive $x \in B$.  Pimsner-Popa index agrees with Jones
index on factors.

If $G$ is an infinite group
with outer action on a $II_1$ factor $X$, then $X^G$ is
not necessarily a factor.  However the Pimsner-Popa index
$[X:X^G]$ is always infinite.
This is trivial if the center of $X^G$ is infinite-dimensional or $\mathbb{C}$,
and may be shown in other cases by examining the inclusions
$p_iX^G \subset p_iXp_i$ for $\{p_i\}$ the minimal central projections
of $X^G$.  It follows that
$M_0^{\omega} \cap M_1' \subset (M_1)_{\omega}$ has infinite 
index if $M_0 \subset M_1=R$ has infinite depth.

\subsection{Bisch-Haagerup subfactors}

Let $H$ and $K$ be finite groups with outer actions on $R$,
and $M_0 \subset M_1$ be the subfactor $R^H \subset R \rtimes K$.
Let the action of $K$ on $R$ be implemented by unitaries
$\{u_k\}$ in the crossed product.
Such group-type subfactors were described by Bisch and Haagerup in \cite{BH}.
We now construct the central sequence subalgebra.
This argument was suggested to the author by V. Jones.

First we note that for a group $X$ with an outer action $\alpha$ on $R$,
$R^{\omega} \cap (R \rtimes X)' = R^{\omega} \cap \{R, \{u_x\}\}'
=R_{\omega} \cap \{u_x\}'$.  The adjoint action of the constant sequences
$(u_x)$ on elements of $R_{\omega}$ is just the pointwise
action of $\alpha_x$ described in the previous section, so this
is $R_{\omega}^X$.  Also every element of
$(R^X)^{\omega} \cap R'$ commutes with the constant sequences
coming from the $u_x$'s, so $(R^X)^{\omega} \cap R'$
is contained in $(R^X)^{\omega} \cap (R \rtimes X)' \subset R_{\omega}^X$
as well.

It follows that $(R^H)^{\omega} \cap (R \rtimes K)'
= ((R^H)^{\omega} \cap R') \cap (R^{\omega} \cap (R \rtimes K)')
\subset R_{\omega}^H \cap R_{\omega}^K$.
Let $G$ now be the image of the free product $H \ast K$ in
${\rm Out} R$, i.e. $G= H \ast K / ({\rm Int} R \cap H \ast K)$.
Then similarly to the previous section we have
$(R^H)^{\omega} \cap (R \rtimes K)' \subset R_{\omega}^G$.

From \cite{BH}, $G$ is infinite if and only if the subfactor
is infinite depth.  In such cases,
$(R^H)^{\omega} \cap (R \rtimes K)'$ is infinite index in $R_{\omega}$,
and hence in $R_{\omega}^K$ as well, since $K$
is finite.  $R_{\omega}^K = R^{\omega} \cap (R \rtimes K)'
\subset (R \rtimes K)_{\omega}$, so if $M^H \subset M\rtimes K$ is infinite
depth then the associated central sequence subalgebra is infinite index.

\section{Towers of finite-dimensional algebras}

\subsection{The canonical shift}

We will now describe the central sequence subalgebra induced by a
commuting-square subfactor.  We begin
with a discussion of towers of finite-dimensional algebras.
This is taken from \cite{JS} and \cite{GHJ}.

Let $A_0 \subset A_1$ be a unital inclusion of finite-dimensional
Von Neumann algebras, with positive trace $\tr$.  We may define the Hilbert
space $L^2(A_1)$ using the inner product $<x,y> = \tr(y^*x)$.
Then the conditional expectation $E_{A_{0}}$ acts on this Hilbert
space.  $A_1$ also acts on $L^2(A_1)$, via left multiplication.
Let $A_2$ be the Von Neumann algebra on $L^2(A_1)$ generated
by $A_1$ and $E_{A_{0}}$.
This is the basic construction on
the inclusion $A_0 \subset A_1$.

Repeating this process (on $A_1 \subset A_2$, etc.) gives the
tower $A_0 \subset A_1 \subset A_2 \subset A_3 \subset ...$.
If the original inclusion is connected (i.e. $A_0 \cap A_1' = \mathbb{C}$)
and the trace on $A_1$ is the unique Markov trace, then there is
a positive trace on the entire tower.  In this case
we may apply the GNS construction to $\cup_i A_i$, and we obtain the
hyperfinite $II_1$ factor $A_{\infty}$.  We will label the
Jones projections $\{e_i\}$ for this tower
according to the convention
$A_i = \{A_{i-1}, e_i\}''$.

From \cite{JS}, these projections have essentially the same properties as the Jones
projections for the tower of factors of~\cite{JoI}.  We mention some of these properties which we will use
in this section.
\begin{itemize}
\item
each $e_i$ has the same trace $\tau$
\item
$e_i e_{i \pm 1} e_i = \tau e_i$
\item
for $x \in A_{i+1}$, $e_{i+2} x e_{i+2} = e_{i+2}E_{A_i}(x)$
\item
for $x \in A_i$, $e_{i+2}x=0$ implies $x=0$
\item
$A_i = A_{i-1}e_iA_{i-1}$
\end{itemize}

From \cite{JS}, if $0 \le i \le j$ then there is an isomorphism from
$A_i' \cap A_j$ to $A_{i+2}' \cap A_{j+2}$.  We now describe this
isomorphism explicitly in terms of the Jones projections.

\begin{lem}

Let $A_0 \subset A_1 \subset A_2 \subset ... \subset A_{\infty}$
be the tower arising from the connected, Markov inclusion
$A_0 \subset A_1$, with Jones projections $\{e_i\}$ labeled
as above.  Let
$w_{ij} = \tau^{-\frac{j-i}{2}}e_{i+2}e_{i+3}e_{i+4}...e_{j+1}e_{j+2}$,
where $\tau$ is the trace of the Jones projections.
Then for all $x \in A_i' \cap A_j$, there is a unique $y \in A_{i+2}' \cap
A_{j+2}$ such that $e_{i+2}y = w_{ij} x w_{ij}^*$, and the map
$\theta_{i,j} : A_i' \cap A_j \rightarrow
A_{i+2}' \cap A_{j+2}$ defined
by $\theta_{i,j}(x)=y$ is an isomorphism.
\end{lem}
\begin{proof}

Let $y$ be an element of $A_{i+2}' \cap A_{j+2}$.
We first note that $$e_{j+2}A_{j+2}e_{j+2} =
e_{j+2} A_{j+1} e_{j+2} A_{j+1} e_{j+2} = e_{j+2} A_j e_{j+2}A_j e_{j+2}
= e_{j+2}A_j$$
Therefore
there is some $x \in A_j$ such that
$e_{j+2}x = w_{ij}^* y w_{ij} \in
e_{j+2} A_{j+2} e_{j+2}$; from
the above properties of the Jones projections, this $x$ is unique.
Define $\rho: A_{i+2}' \cap A_{j+2} \rightarrow A_j$
by $e_{j+2}\rho(y) = w_{ij}^* y w_{ij}$.  We will show that this
map is a $^*$-isomorphism into $A_i' \cap A_j$, and compute
$\theta_{i,j}$ as its inverse.

It follows from the definition of $\rho$ that
$\rho(y) = \tau^{-1}E_{A_j}(w_{ij}^* y w_{ij})$.  Therefore $\rho$ respects
the $^*$ operation.  Furthermore, since $A_i \subset A_j$,
and both $w_{ij}$ and $y$ commute with $A_i$, $\rho(y)$ commutes
with $A_i$ for all $y$ and $\rho(A_{i+2}' \cap A_{j+2}) \subset
A_i' \cap A_j$.

We may compute as well
$$e_{j+2} \rho(y_1) \rho(y_2) = e_{j+2} \rho(y_1) e_{j+2} \rho(y_2)
= w_{ij}^* y_1 w_{ij} w_{ij}^* y_2 w_{ij}$$
But $w_{ij} w_{ij}^* = e_{i+2}$ by the properties
of the Jones projections. $e_{i+2}$ commutes with $y_1$,
and $w_{ij}^* e_{i+2} = w_{ij}^*$, so this is $w_{ij}^* y_1y_2 w_{ij} =
e_{j+2} \rho(y_1y_2)$.
Therefore $\rho$ is a homomorphism.

Now let $y$ be a positive element of $A_{i+2}' \cap A_{j+2}$.
$\rho(y) = 0$ only if $w_{ij}^* y w_{ij} = 0$, and
$$\tr(w_{ij}^* y w_{ij}) =
\tr(y w_{ij} w_{ij}^*) = \tr(y e_{i+2})$$
This is equal to
$\tr(E_{A_{i+2}}(y e_{i+2})) = \tr(E_{A_{i+2}}(y)e_{i+2})$.
Since $y>0$ and $y$ commutes with $A_{i+2}$, $E_{A_{i+2}}(y)$
is a positive element of the center $Z(A_{i+2})$.  $e_{i+2}$
has full central support in $A_{i+2}$, so for any $a>0$ in
$Z(A_{i+2})$, $\tr(e_{i+2}a)>0$.
So $\tr(e_{i+2}E_{A_{i+2}}(y))=\tr(e_{i+2}y)=\tr(\rho(y))$ is positive.
This means that for any $z \in A_{i+2}' \cap A_{j+2}$,
$\tr(\rho(z)^*\rho(z)) = \tr(\rho(z^*z))>0$, and $\rho$ is an injective
homomorphism from $A_{i+2}' \cap A_{j+2}$ onto $A_i' \cap A_j$.

From \cite{JS}, the algebras $A_{i+2}' \cap A_{j+2}$ and $A_i' \cap A_j$
are isomorphic, so $\rho$ is an isomorphism.  It follows that
for $x \in A_i' \cap A_j$, $x = \rho(y)$ for some unique $y \in
A_{i+2}' \cap A_{j+2}$.

$x = \rho(y)$ means that $e_{j+2} x = w_{ij}^* y w_{ij}$.  Conjugating
by $w_{ij}$ and its adjoint shows that this is true
if and only if $ e_{i+2}y = w_{ij}x w_{ij}^* $, so there is also
a unique $y$ obeying this second relation.  If we define a map
$\theta_{i,j}$ from $A_i' \cap A_j$
to $A_{i+2}' \cap A_{j+2}$ by $\theta_{i,j}(x)=y$, then $\theta_{i,j}
= \rho^{-1}$ and is an isomorphism as desired.
\end{proof}

This map $\theta$ is the canonical shift on the tower of finite-dimensional
algebras (c.f.~\cite{Oc}).
We may likewise construct the iterated shift $\theta_{i,j}^k$
as the product
$$\theta^k_{i,j}=\theta_{i+2(k-1),j+2(k-1)}\theta_{i+2(k-2),j+2(k-2)}
...\theta_{i+2,j+2}\theta_{i,j}$$
This is an isomorphism from $A_i' \cap A_j$ to $A_{i+2k}' \cap A_{j+2k}$.

The asymptotic properties of the trace on the tower
will be important later on,
so we mention some results derived from Perron-Frobenius theory,
following \cite{JS}.

For a finite-dimensional Von Neumann algebra A with minimal central
projections $\{p_1 ... p_n\}$, we have a trace vector $\mathbf{t}$ and a size
vector $\mathbf{s}$; i.e., the matrix subalgebra $p_i A$ is $s_i$ by $s_i$,
and the trace of a minimal projection in this subalgebra is $t_i$.
From normalization of the trace, we must have $<\mathbf{s}, \mathbf{t}> =
\sum_{i=1}^n s_i t_i=1$.

Let $\mathbf{s}^{(i)}$ and $\mathbf{t}^{(i)}$ be the size and trace vectors for $A_{i}$.
From \cite{JS}, the inclusion matrices of $A_{n} \subset A_{n+1}$
and $A_{n+1} \subset A_{n+2}$ are transposes of each other for all $n$.
Let $\Lambda$ be the inclusion matrix of $A_i \subset A_{i+1}$; then
this means that there is a labeling of the central projections of the
$A_i$'s such that for all $k \ge 0$ we have
$\mathbf{s}^{(i+2k+1)} = \Lambda^T \mathbf{s}^{(i+2k)}$ and
$\mathbf{s}^{(i+2k+2)}= \Lambda \mathbf{s}^{(i+2k+1)}$.
Since $A_i \subset A_{i+1}$ is a connected inclusion, all entries
of $(\Lambda \Lambda^T)^k$ are positive for $k$ sufficiently large,
and this matrix has a unique Perron-Frobenius eigenvector $v$
with $||v||=1$ and $v_i > 0$.
The Perron-Frobenius eigenvalue is $||\Lambda\Lambda^T||=\tau^{-1}=
\tr(e_i)^{-1}$.

We likewise have $\mathbf{t}^{(i+2k)} = \Lambda \mathbf{t}^{(i+2k+1)}$,
$\mathbf{t}^{(i+2k+1)} = \Lambda^T \mathbf{t}^{(i+2k+2)}$. 
From \cite{JS}, the Markov condition implies that
$\mathbf{t}^{(i+2k)}$ is a positive scalar multiple of $v$ for all $k$,
so for all $k \ge 0$ we have $\mathbf{t}^{(i+2k+2)} = \tau \mathbf{t}^{(i+2k)}$.

It then follows from Perron-Frobenius theory that
$\lim_{k \rightarrow \infty} \mathbf{s}^{(i+2k)}\tau^k = \gamma \mathbf{t}^{(i)}$,
for some positive scalar $\gamma$.  The normalization condition and
the above expression for $\mathbf{t}^{(i+2k)}$
implies that $\gamma = ||\mathbf{t}^{(i)}||^{-2}$.  With this result we can
describe the asymptotic behavior of the trace on
$A_{i+2k}' \cap A_{j+2k}$.

\begin{lem}
Let $A_0 \subset A_1 \subset ... \subset A_{\infty}$ be a Jones
tower as above.  Fix $0 \le i \le j$.
There exists $\epsilon>0$ such that for all $k \ge 0$,
and all nonzero projections $p \in A_{i+2k}' \cap A_{j+2k}$, we have
$\tr(p) \ge \epsilon$.
\end{lem}
\begin{proof}
For $k \ge 0$, we take the minimal central projections of $A_{i+2k}$,
$A_{j+2k}$ to be $\{p^k_x\}$, $\{q^j_y\}$.  As above, we may
use the same set of labels for all $k$.
Every minimal central projection $r^k_{xy}$ in $A_{i+2k}' \cap A_{j+2k}$
is of the form $p^k_xq^k_y$.
Let $\tilde{r}^k_{xy}$ be an arbitrary minimal projection in
the matrix algebra
$r^k_{xy}(A_{i+2k}' \cap A_{j+2k})$.
From \cite{JS} we may compute
$\tr(\tilde{r}^k_{xy})=s^{(i+2k)}_x t^{(j+2k)}_y$.

Since $\mathbf{t}^{(j+2k)} = \tau^k \mathbf{t}^{(j)}$ and $\mathbf{s}^{(i+2k)} \tau^k$ approaches
$\mathbf{t}^{(i)} ||\mathbf{t}^{(i)}||^{-2}$, this trace approaches the limit $t^{(i)}_xt^{(j)}_y||\mathbf{t}^{(i)}||^{-2}$ as
$k$ goes to $\infty$.
Then for fixed $i$ and $j$, the sequence $(\tr(\tilde{r}^k_{xy}))$
consists of positive numbers approaching a positive
value, and is therefore bounded away from zero.  Since the centers
$Z(A_i)$ and
$Z(A_j)$ are finite-dimensional, there is some number $\epsilon>0$
which bounds the traces of these minimal projections
away from zero for all $x$, $y$.
We then have $\epsilon \le \tr(\tilde{r}^k_{xy})$ for all $k$, $x$, $y$.
Since the $\tilde{r}^k_{xy}$'s are minimal projections in every
central component of $A_{i+2k}' \cap A_{j+2k}$, it follows
that $\epsilon$ is less than or equal to the trace of any nonzero projection
in any of these algebras.
\end{proof}

\subsection{The grid of finite-dimensional algebras}

Let $$\commsq{A_{00}}{A_{10}}{A_{01}}{A_{11}}$$ be a quadrilateral
of finite-dimensional Von Neumann algebras, with trace.
We may construct the Hilbert space
$L^2(A_{11})$ as in the previous section.
This quadrilateral is a commuting square if
$E_{A_{10}}$ commutes with $E_{A_{01}}$ on this Hilbert space.

In \cite{GHJ} (see also \cite{JS}),
the authors construct a subfactor from such
a commuting square.  The commuting square must have some additional
properties: every constituent inclusion must be connected
and have the unique Markov trace.  Additionally, the square must
be symmetric.
Several equivalent conditions for symmetry are discussed in \cite{JS};
one of them is $A_{10}A_{01}=A_{11}$, which we will take as our
definition.

All of our commuting squares in this paper will be connected,
symmetric, and Markov.
The following construction is taken from \cite{JS}.

Starting with such a commuting square, we may build the Jones
tower $A_{01} \subset A_{11} \subset A_{21} \subset ... \subset A_{\infty 1} = M_1$,
with Jones projections $\{e_i\}$ labeled as in the previous section.
Then we may define $A_{i0}$ inductively for $i \ge 2$ by
$A_{i0} = \{A_{i-1,0}, e_i\}''$.  From \cite{JS}, the
$A_{i0}$'s form a Jones tower as well, and
$M_0=A_{\infty 0} = \overline{\cup_i A_{0i}}^{st}$
is a subfactor of the hyperfinite
$II_1$ factor $M_1$.

The index $[M_1:M_0]$ will be finite, and we may construct
the tower of $II_1$ factors $M_0 \subset M_1 \subset M_2 \subset ...$.
We will label the vertical Jones projections for this tower by
$M_j = \{M_{j-1}, f_j\}''$, $j \ge 2$.

We may then construct a grid of finite-dimensional algebras,
defining $A_{ij} \subset M_j$ by $\{A_{i,j-1},f_j\}''$. 
$A_{0j} \subset A_{1j} \subset A_{2j} \subset ... \subset A_{\infty j} = M_j$
is again a Jones tower, and we may use the same Jones projections
$\{e_i\}$ for any value of $j$.

This implies that
for $0 \le j \le k$, $M_j \subset M_k$ is also a commuting-square
subfactor.  It may be shown that the quadrilateral
$$\commsq{A_{0j}}{A_{1j}}{A_{0k}}{A_{1k}}$$ commutes and is
symmetric (i.e., $A_{0k}A_{1j}=A_{1k}$) by
induction on $k-j$.
The inclusion $A_{0k} \subset A_{1k}$ is Markov, so the entire square
is Markov \cite{JS}, and connectedness follows from connectedness
of the original commuting square.  So this quadrilateral generates
a commuting-square subfactor.  Since the Jones projections for the tower of
 $A_{0k} \subset A_{1k}$ are the $e_i$'s, this subfactor may
 be taken to be $(A_{1j} \cup \{e_i\})'' \subset
 (A_{1k} \cup \{ e_i\})'' = M_j \subset M_k$ .

We may likewise consider the vertical tower
$A_{0i} \subset A_{1i} \subset ...$.  From \cite{JS} this is again
a Jones tower with the same Jones projections $\{f_j\}$ for any
value of $i$.  The vertical limits
$P_i=A_{\infty i} = \overline{\cup_j A_{ji} }^{st}$ are
$II_1$ factors, with $P_0 \subset P_1 \subset P_2 \subset ...$
a Jones tower.

This grid may be used to explicitly compute the higher relative
commutants of a commuting-square subfactor.  Ocneanu's
compactness argument, described in \cite{JS}, shows that
$M_0' \cap M_k = A_{10}' \cap A_{0k}$.  An analogous statement
may be made about the vertical relative commutants
$P_0' \cap P_k$.

It follows from Perron-Frobenius arguments similar to those used in Lemma
2 that all commuting-square subfactors are extremal.  This result will
be necessary
for demonstrating the desired inequalities on index in section 4.

\begin{lem}
Let the grid of algebras $\{ A_{ij} \}$ be as above.
Let $x$ be an element of $A_{0k}$.
Then for $j \le k$,
$\lim_{i \rightarrow \infty} E_{A_{ij}' \cap A_{ik}}(x)$ exists
and equals $E_{M_j' \cap M_k}(x)$.
\end{lem}
\begin{proof}

Fix $x \in A_{0k}$.

We consider the projections $E_{A_{ij}' \cap M_k}$, acting
on $L^2(M_k)$.  This is a decreasing series of orthogonal projections
onto closed subspaces of $L^2(M_k)$, and therefore strongly approaches
the orthogonal projection onto the intersection of these subspaces.
So $\lim_{i \rightarrow \infty} E_{A_{ij}' \cap M_k}(x)$
exists and equals $E_{\cap_i(A_{ij}' \cap M_k)}(x)$.
We may compute
$\cap_i (A_{ij}' \cap M_k) = (\cup_i A_{ij})' \cap M_k
= (\overline{(\cup_i A_{ij})}^{st})' \cap M_k = M_j' \cap M_k$.

The quadrilateral
$$\commsq{A_{ij}' \cap A_{ik}}{A_{ij}' \cap M_k}
{A_{ik}}{M_k}$$
commutes, since $E_{A_{ik}}(A_{ij}' \cap M_k)$ commutes
with $A_{ij}$.  Therefore
$E_{A_{ij}' \cap M_k }(x) = E_{A_{ij}' \cap A_{ik}}(x)$,
since $x \in A_{0k} \subset A_{ik}$.
It follows that $\lim_{i \rightarrow \infty}
E_{A_{ij}' \cap A_{ik}}(x)$ exists and equals $E_{M_j' \cap M_k}$,
as desired.
\end{proof}

\begin{thm}
All commuting-square subfactors are extremal.
\end{thm}
\begin{proof}

Let the grid of algebras $\{A_{ij} \}$ be as above.  Let
$\tau_h = \tr(e_i)$ and $\tau_v = \tr(f_i)$ be the Markov constants
for the horizontal and vertical inclusions respectively, and
let $A_{ij}$ have size vector $\mathbf{s}^{(ij)}$ and trace vector $\mathbf{t}^{(ij)}$.
From \cite{JS}
we may label the central projections of the
$A_{ij}'s$ so that $\mathbf{t}^{(i+2k,j+2l)} = \mathbf{t}^{(ij)} \tau_v^k \tau_h^l$.
From the previous section, with this labeling we have
$$\lim_{k \rightarrow \infty} \mathbf{s}^{(i+2k,j+2l)} \tau_h^k =
\mathbf{t}^{(i,j+2l)} ||\mathbf{t}^{(i,j+2l)}||^{-2} =
\mathbf{t}^{(ij)}||\mathbf{t}^{(ij)}||^{-2} \tau_v^{-l}$$

Now we consider the inclusion of algebras $A_{2k,0} \subset
A_{2k,1} \subset A_{2k,2}$.  Let the minimal
central
projections of these three algebras be
$\{p^{(k)}_x\}$, $\{q^{(k)}_y\}$, and $\{\tilde{p}^{(k)}_x\}$
respectively, labeled as above.  From \cite{JS},
we then have $p^{(k)}_x f_2 = \tilde{p}^{(k)}_x f_2$.
The minimal central projections of $A_{2k,0}' \cap A_{2k,1}$ are all
of the form $p^{(k)}_x q^{(k)}_y$, and those of
$A_{2k,1}' \cap A_{2k,2}$ are of the form
$\tilde{p}^{(k)}_x q^{(k)}_y$.  Our choice of labels lets
us use the same indices $(x,y)$ for all $k$.

Let $r^{(k)}_{xy}$
be any
minimal projection in $p^{(k)}_x q^{(k)}_y(A_{2k,0}' \cap A_{2k,1})$.
Then $\tr(r^{(k)}_{xy}) = s^{(2k,0)}_x t^{(2k,1)}_y$.
This means that
$\lim_{k \rightarrow \infty} \tr(r^{(k)}_{xy}) =
t^{(00)}_x t^{(01)}_y ||\mathbf{t}^{(01)}||^{-2}$.
Similarly,
if $\tilde{r}^{(k)}_{xy}$ is an arbitrary minimal projection in
$\tilde{p}^{(k)}_x q^{(k)}_y(A_{2k,1}' \cap A_{2k,2})$,
then
$\lim_{k \rightarrow \infty} \tr(\tilde{r}^{(k)}_{xy}) =
t^{(00)}_x t^{(01)}_y ||\mathbf{t}^{(00)}||^{-2} \tau_v$.

If $\Lambda_v$ is the inclusion matrix for
$A_{00} \subset A_{01}$, then
$$||\mathbf{t}^{(00)}||^{-2} = ||\Lambda_v^T \mathbf{t}^{(01)}||^{-2} =
<\Lambda_v \Lambda_v^T \mathbf{t}^{(01)}, \mathbf{t}^{(01)}>^{-1} =
||\mathbf{t}^{(01)}||^{-2}\tau_v^{-1}$$
So these two limits are the same, and
for all $x, y$ and $\epsilon>0$
there exists $k_0$ such that
$k>k_0$ implies
$||\tr(\tilde{r}^{(k)}_{xy})-\tr(r^{(k)}_{xy})||<\epsilon$.

From \cite{JS}, $A_{2k,2}$ acts on the Hilbert
space $L^2(A_{2k,1})$, where $A_{2k,1}$ acts
by left multiplication and $f_2$ is the conditional expectation onto
$A_{2k,0}$.  If $J$ is the order-2 anti-linear anti-isometry
given by $J(x)=x^*$, then on this Hilbert space
$JA_{2k,1}J = A_{2k,1}'$, $JA_{2k,0}'J = A_{2k,2}$.  It
follows that $J(A_{2k,0}' \cap A_{2k,1})J = A_{2k,1}' \cap A_{2k,2}$.
If we define $\phi(x)$ on $A_{2k,0}' \cap A_{2k,1}$ by
$\phi(x) = J x^* J$, then $\phi$ is a $^*$-isomorphism from
$A_{2k,0}' \cap A_{2k,1}$ onto $A_{2k,1}' \cap A_{2k,2}$.

With this labeling we have
$\phi(p_x q_y) =
\phi(p_x) \phi(q_y) = \tilde{p}_x q_y$ \cite{JS}.  Therefore $\phi$
sends a minimal projection in $p_x q_y (A_{2k,0}' \cap A_{2k,1})$
to one in $\tilde{p}_x q_y (A_{2k,1}' \cap A_{2k,2})$.

For any $\epsilon>0$,
the above argument on convergence 
of traces (and Lemma 2) gives $k_0$ such that
$k \ge k_0$ implies $|\tr(\phi(r)) - \tr(r)| < \epsilon||r||_2$
for all minimal projections $r$ in
$A_{2k,0}' \cap A_{2k,1}$.
Since
$A_{2k,0}' \cap A_{2k,1}$ is spanned by minimal
projections, and is of fixed finite dimension, for $k$ sufficiently
large
we then have
$|\tr(\phi(v))-\tr(v)| <\epsilon||v||_2$ for all
$v \in A_{2k,0}' \cap A_{2k,1}$.

Now choose $x \in A_{2k,0}' \cap A_{2k,1}$, and let
$z$ be an arbitrary element of $L^2(A_{2k,1})$.
We compute $(x f_2)(z) = x E_{A_{2k,0}}(z)$. 
$x$ commutes with $A_{2k,0}$, so this is equal to
$$E_{A_{2k,0}}(z) x = (x^* E_{A_{2k,0}}(z)^*)^* = (Jx^*Jf_2)(z) =
(\phi(x)f_2)(z)$$
 Since $\phi(x) f_2$ and $x f_2$ agree on
$L^2(A_{2k,1})$, they are the same element of $A_{2k,0}' \cap A_{2k,2}$.

This means that $\tr(\phi(x) f_2) = \tr(x f_2)$, which is equal
to $\tr(x) \tau_v$ by the Markov property of the trace.
For $k$ sufficiently large, this is within $\epsilon||x||_2$ of
$\tr( \phi(x) ) \tau_v = \tr(\phi(x) (\tau_v 1))$, for any
$\phi(x) \in A_{2k,1}' \cap A_{2k,2}$.

This implies
$||E_{A_{2k,1}' \cap A_{2k,2}}(f_2)-\tau_v 1||_2 < \epsilon$
for $k$ sufficiently large, so $\lim_{k \rightarrow \infty}
E_{A_{2k,1}' \cap A_{2k,2}}(f_2) = \tau_v 1$.
From Lemma 3, this means that $E_{M_1' \cap M_2}(f_2) = \tau_v 1$,
and so $M_0 \subset M_1$ is extremal.
\end{proof}

\section{Commuting-square subfactors and central sequences}

\subsection{Remarks on index}

The iterated canonical shift from the previous section allows
us to give a bound for the norm of certain central sequence subalgebras.

\begin{lem}
Let $M_0 \subset M_1$ be the subfactor obtained by iterating
the basic construction on the commuting square
$$\commsq{A_{00}}{A_{10}}{A_{01}}{A_{11}}$$
Let the grid of algebras be
$\{A_{ij}\}$ as in section 3.2.  Let $p$
be a projection in $A_{0k}' \cap A_{ik}$.
Then $[(M_k)_{\omega}:M_0^{\omega} \cap M_1'] \ge
||E_{A_{00}' \cap A_{i0}}(p)||^{-1}$.
\end{lem}
\begin{proof}

The conditions of Lemma 1 are satisfied by the Jones tower
$A_{0k} \subset A_{1k} \subset ... \subset M_k$.
So there are isomorphisms $\{\theta^l_{i,j}\}$ from
$A_{ik}' \cap A_{jk}$ to $A_{i+2l,k}' \cap A_{j+2l,k}$.
Isomorphisms $\{\psi^l_{i,j}\}$ likewise exist
from $A_{i0}' \cap A_{j0}$ to $A_{i+2l,0}' \cap A_{j+2l,0}$,
as in the lemma.

For $p$ a projection in $A_{0k}' \cap A_{ik}$,
let $\tilde{p}$ be the element of $M_k^{\omega}$ given
by $\tilde{p}_l = \theta^l_{0,j}(p)$.  This is a projection.

By construction $\tilde{p}$ commutes with $A_{2l,k}$ for
all $l$.  Since these algebras generate $M_k$, $\tilde{p}$
is contained in $(M_k)_{\omega}$.  Since 
$\theta^l_{0,j}(p) \in A_{2l,k}' \cap A_{j+2l,k}$,
from Lemma 2 there
exists $\epsilon>0$ such that
$\tr(\theta^l_{0,j})(p) \ge \epsilon$ for all $l$.  Therefore
$\tilde{p} \neq 0$ as an element of $(M_k)_{\omega}$.

Let $x$ be an element of $A_{i,k}' \cap A_{j,k}$.
From section 3.1, $\theta^1_{i,j}(x)=\theta_{i,j}(x)$
is defined as the unique element
of $A_{i+2,k}' \cap A_{j+2,k}$ such that
$e_{i+2} \theta_{i,j}(x) = w_{ij} x w_{ij}^*$.
$w_{ij}$
is a word in the horizontal
Jones projections, which are contained in $M_0$.
So applying $E_{M_0}$ to both sides, we find that
$e_{i+2} E_{M_0}(\theta_{i,j}(x)) = w_{ij} E_{M_0}(x) w_{ij}^*$.

For any $n \ge 0$, it may be shown by induction on $k$ that
$E_{M_0}(A_{kn})=A_{0n}$.
Therefore
$E_{M_0}(\theta_{i,j}(x)) \in A_{i+2,0}' \cap A_{j+2,0}$
and $E_{M_0}(x) \in A_{i,0}' \cap A_{j,0}$.
But for $y \in A_{i,0}' \cap A_{j,0}$,
$\psi_{i,j}(y)$ is defined as the unique
element of $A_{i+2,0}' \cap A_{j+2,0}$
obeying the relation $e_{i+2}\psi_{i,j}(y) = w_{ij} y w_{ij}^*$.
Therefore for all $i \le j$, and all $x \in A_{ik}' \cap A_{jk}$,
$\psi_{i,j}(E_{M_0}(x)) = E_{M_0}(\theta_{i,j}(x))$.

From the definitions of the composite operators
$\psi^l_{0,j}$ and $\theta^l_{0,j}$ in section 3.1,
it follows that $E_{M_0} (\theta^l_{0,j} (p)) =
\psi^l_{0,j} (E_{M_0}(p))$.

We may compute the conditional expectation onto
$M_0^{\omega}$ by applying $E_{M_0}$ pointwise.
Therefore $||E_{M_0^{\omega}}(\tilde{p})_l ||=
||E_{M_0}(\Theta_{0,j}^l(p))|| =||\psi_{0,j}^l(E_{M_0}(p))||$.  Since
isomorphisms preserve $\infty$-norm, this is
$||E_{M_0}(p)||$.

All components of $E_{M_0^{\omega}}(\tilde{p})$ have
norm at most $||E_{M_0}(p)||$, so the same is true
of $E_{M_0^{\omega}}(\tilde{p})$ itself.  Since
$M_0^{\omega} \cap M_1'$ is a Von Neumann subalgebra of
$M_0^{\omega}$, we have
$||E_{M_0^{\omega} \cap M_1'}(\tilde{p})|| \le ||E_{M_0}(p)||$
as well.  From \cite{PP}, since $\tilde{p}$ is a nonzero projection
in $(M_1)_{\omega}$, this gives
$$[(M_1)_{\omega}:M_0^{\omega}\cap M_1'] \ge ||E_{M_0}(p)||^{-1}
= ||E_{A_{00}' \cap A_{i0}}(p)||^{-1}$$
as desired.
\end{proof}

The above argument only bounds $[(M_k)_{\omega}:M_0^{\omega} \cap M_1']$,
which is not quite what we want.  Some additional lemmas on index
will allow us to show that this bound also applies to
$[(M_1)_{\omega}:M_0^{\omega} \cap M_1']$.

\begin{lem}
Let $X \subset Y$ be a $II_1$ subfactor, with $X$ hyperfinite.
Then $X' \cap Y^{\omega}$ is a $II_1$ factor.
\end{lem}
\begin{proof}
For $L \subset P \subset Q$ $II_1$ factors, and $L$ hyperfinite,
the central freedom lemma (see \cite{KaC}) states that
$(L' \cap P^{\omega})' \cap Q^{\omega} = L \vee (P' \cap Q)^{\omega}$.

We apply this lemma to the inclusion $X \subset Y \subset Y$,
obtaining $(X' \cap Y^{\omega})' \cap Y^{\omega} = X \vee (Y' \cap Y)^{\omega}$.
$Y$ is a factor, so this is just $X$.  This means
that $(X' \cap Y^{\omega})' \cap Y^{\omega} \cap X' = X \cap X' = \mathbb{C}1$,
since $X$ is a factor as well.  So
$(X' \cap Y^{\omega})$ has trivial center, and is a factor.

$X' \cap Y^{\omega}$ is contained in the $II_1$ factor $Y^{\omega}$,
and contains the $II_1$ factor $X_{\omega}$.  It is therefore
of type $II_1$.
\end{proof}

\begin{lem}
Let $Y \subset Z$ be a finite-index $II_1$ subfactor.  Let
$X$ be a Von Neumann subalgebra of $Y$, with
finite Pimsner-Popa index $[Y:X]$.  Then
$[Z:X]=[Z:Y][Y:X]$.
\end{lem}
\begin{proof}
By applying the downward basic construction (c.f. \cite{JoI})
we may
find a projection $e \subset Z$ such that $E_Y(e) = [Z:Y]^{-1}1$
and $Y_0 = e' \cap Y$ is a $II_1$ factor.
Also, from \cite{PP}, for all $\epsilon$
we have a positive element $q \in Y$ such that
$||E_X(q)||$ is not greater than $([Y:X]^{-1} + \epsilon)q$.

$Y_0 \subset Y_1$ is a finite-index
$II_1$ subfactor, and $q$ may be approximated
to arbitrary precision in 2-norm by a finite linear combination
of orthogonal projections.  Therefore there
exists a unitary $u \in Y$ such that
$uqu^* \in Y_0$.  Since $uqu^*$ commutes with $e$, it follows
that $u^*eu$ commutes with $q$.  This means that
$p = u^* e u q $ is a positive element of $Z$.

Now we compute $||E_X(p)||$.  This is equal to $||E_X(E_Y(p))||$.
$E_Y(p) = E_Y(u^* e u q) = u^* E_Y(e) u q = [Z:Y]^{-1}q$.
$E_X([Z:Y]^{-1}q) = [Z:Y]^{-1}E_X(q)$, which is not
greater than $[Z:Y]^{-1}( [Y:X]^{-1} + \epsilon)q$.
Therefore $[Z:X] \ge [Z:Y][Y:X]$.

For any positive element $a \in Z$, we have
$E_X(a) = E_X(E_Y(a))$.  By definition of Pimsner-Popa
index, $E_Y(a) \ge [Z:Y]^{-1}a$.
So
$$E_X(E_Y(a)) \ge [Y:X]^{-1}E_Y(a) \ge [Z:Y]^{-1}[Y:X]^{-1}a$$
Since this is true for all $a \ge 0$ in $Z$, we have
$[Z:X] \le [Z:Y][Y:X]$ as well.  This gives the desired equality.
\end{proof}

\begin{lem}
Let $$\commsq{A}{B}{C}{D}$$ be a quadrilateral of Von Neumann
algebras.  Let $B$, $C$, and $D$ be $II_1$ factors,
with $[D:B] \le [D:C]$ and $[D:B] < \infty$.  Then
$[C:A] \le [B:A]$.
\end{lem}
\begin{proof}
First we consider the case $[D:A] = \infty$.  From [PP],
if $X \subset Y \subset Z$ are Von Neumann algebras with finite
trace,
and $[Z:X] = \infty$, then either $[Y:X]$ or $[Z:Y]$ must
be infinite.  Since $[D:B] < \infty$, we must have
$[B:A] = \infty$, giving the desired inequality for any value
of $[C:A]$.

Now let $[D:A]$ be finite.  This implies that all four inclusions
of the above quadrilateral are finite index.  From Lemma 6,
since $D$ and $B$ are factors we have $[D:A] = [D:B][B:A]$.
Likewise $[D:A] = [D:C][C:A]$.  Since all indices are positive
real numbers here, we may compute
$[C:A] = \frac{[D:B][B:A]}{[D:C]}$.  $[D:B] \le [D:C]$ by
hypothesis, so $[C:A] \le [B:A]$ in this case as well.
\end{proof}

With these results, we can start to approximate the indices of various kinds
of central sequence inclusions.

\begin{lem}
Let $M_0 \subset M_1$ be a finite-index $II_1$ factor,
with $X$ a hyperfinite subfactor of $M_0$.  Then
$[X' \cap M_1^{\omega}:X' \cap M_0^{\omega}]=[M_1:M_0]$.
\end{lem}
\begin{proof}
We first consider the quadrilateral of Von Neumann algebras
$$\commsq{X' \cap M_0^{\omega}}{X' \cap M_1^{\omega}}
{M_0^{\omega}}{M_1^{\omega}}$$
This is a commuting square since $X \subset M_0^{\omega}$.
So for any $a>0$ contained in
$X' \cap M_1^{\omega}$
we have $E_{X' \cap M_0^{\omega}}(a) = E_{M_0^{\omega}}(a)
\ge [M_1^{\omega}:M_0^{\omega}]^{-1}a$.  This
means that $[X' \cap M_1:X' \cap M_0] \ge [M_1^{\omega}:M_0^{\omega}]$,
which is equal to $[M_1:M_0]$ from \cite{PP}.

From the downward basic construction there exists $e \in M_1$ with
$E_{M_0}(e) = [M_1:M_0]^{-1}$, $e' \cap M_0 = M_{-1}$, where
$M_{-1}$ is a $II_1$ factor.

Let $X = \overline { \{\cup_i X_i\} }^{st}$, where each
$X_i$ is a matrix algebra and $X_i \subset X_{i+1}$.
There is a matrix algebra $A_i \subset M_{-1}$ of the same size as
$X_i$.  Any two matrix algebras of the same size in a $II_1$
factor are unitarily equivalent, so there is a unitary
$u_i \in M_0$ with $u_i X_i u_i^* \subset M_{-1}$.

Since $u_i X_i u_i^*$ commutes with $e$, it follows that
$u_i^* e u_i$ commutes with $X_i$.  The sequence
$\tilde{e}$ defined by
$\tilde{e}_i = u_i^* e u_i$ gives a projection in $M_1^{\omega}$.
This sequence asymptotically commutes with every $X_i$.
Since the $X_i$'s are dense in $X$, $\tilde{e} \in X' \cap M_1^{\omega}$.

As above,
$E_{X' \cap M_0^{\omega}}(\tilde{e}) = E_{M_0^{\omega}}(\tilde{e})$.
We may compute this conditional expectation by
applying $E_{M_0}$ pointwise, obtaining
$(E_{M_0^{\omega}}(\tilde{e}))_i =
E_{M_0}(u_i^* e u_i) = u_i^* E_{M_0}(e) u_i = [M_1:M_0]^{-1}1$.
This is a constant sequence, and so $
E_{X' \cap M_0^{\omega}}(\tilde{e}) = [M_1:M_0]^{-1}1$ as an element
of $X' \cap M_0^{\omega}$.  It follows that
$[X' \cap M_1^{\omega}:X' \cap M_0^{\omega}] \ge [M_1:M_0]$.
The reverse inequality has been shown above, so the two indices
are equal.
\end{proof}

\begin{lem}
Let $M_0 \subset M_1$ be an extremal finite-index hyperfinite
$II_1$ subfactor.  Then
$[M_0' \cap M_1^{\omega}:(M_1)_{\omega}] \ge [M_1:M_0]$.
\end{lem}
\begin{proof}
Let $X \subset Y \subset Z$ be an inclusion
of $II_1$ factors, with $X$ hyperfinite.
We consider the quadrilateral of Von Neumann algebras
$$\commsq{Y' \cap Z}{X' \cap Z}
{Y' \cap Z^{\omega}}{X' \cap Z^{\omega}}$$
$X' \cap Z^{\omega}$ is a $II_1$ factor by Lemma 5, so
there is a unique trace on this quadrilateral, and we may
compute conditional expectations.
Since $Y \subset Z$, $E_Z(Y' \cap Z^{\omega})$ is contained
in $Y' \cap Z$, and so $E_{X' \cap Z}(Y' \cap Z^{\omega})$
is as well.  Therefore the quadrilateral is a commuting square, and
in general for $a \in X' \cap Z$,
 $E_{Y' \cap Z^{\omega}}(a) = E_{Y' \cap Z}(a)$.

Let $M_2$ be obtained by applying the basic construction to
$M_0 \subset M_1$, with Jones projection $e \in M_2$.
We consider the quadrilateral of $II_1$ factors
$$\commsq{(M_1)_{\omega}}{M_0' \cap M_1^{\omega}}
{M_1' \cap M_2^{\omega}}{M_0' \cap M_2^{\omega}}$$

We may embed $e$ in $M_2^{\omega} \cap M_0'$ as a constant sequence.
The above argument shows that
$E_{M_1' \cap M_2^{\omega}}(e) = E_{M_1' \cap M_2}(e)$.  This
is $[M_2:M_1]^{-1}1$
since $M_0 \subset M_1$ is extremal by hypothesis.  Therefore
$[M_0' \cap M_2^{\omega}:M_1' \cap M_2^{\omega}] \ge [M_2:M_1]$.

From Lemma 8, we have
$[M_0' \cap M_1^{\omega}:M_0' \cap M_2^{\omega}]=[M_2:M_1]$.
$[M_2:M_1] = [M_1:M_0]$ \cite{JoI}, which is finite by hypothesis.
So all the assumptions of Lemma 7 are satisfied,
implying that $[M_0' \cap M_1^{\omega}:(M_1)_{\omega}] \ge
[M_1' \cap M_2^{\omega}:(M_1)_{\omega}]$.
Again by Lemma 8, 
$[M_1' \cap M_2^{\omega}:(M_1)_{\omega}] = [M_2:M_1]$, giving
the desired inequality.
\end{proof}

\subsection{Infinite-depth central sequence subfactors}

We have not yet used the fact that our initial subfactor is of infinite
depth.  A result of Sato allows us to use this
to describe the asymptotic
behavior of $[(M_k)_{\omega}:M_0' \cap M_1^{\omega}]$.

\begin{lem}
Let $M_0 \subset M_1$ be a commuting-square subfactor of infinite
depth, with Jones tower $M_0 \subset M_1 \subset M_2 \subset ...$.
Then
$\lim_{k \rightarrow \infty} [(M_k)_{\omega}:M_0' \cap M_1^{\omega}]= \infty$.
\end{lem}
\begin{proof}
Let $M_0 \subset M_1$ be generated by the commuting square
$$\commsq{A_{00}}{A_{10}}{A_{01}}{A_{11}}$$
with the grid
of algebras $\{A_{ij}\}$
as in section 3.2.  The centers
of the $A_{ij}$'s have bounded dimension; let this bound be $L$.
Then ${\rm dim} Z(A_{ij}' \cap A_{kl})<L^2$ for any
$0 \le i \le k, 0 \le j \le l$.

Choose $\epsilon>0$.

Since $M_0 \subset M_1$ is of infinite depth,
the vertical subfactor
$A_{0\infty} \subset A_{1\infty} = P_0 \subset P_1$ is also
of infinite depth \cite{Sa}.  Therefore the central
dimension of $P_0' \cap P_i$ increases
without limit as $i$ goes to infinity.
By Ocneanu compactness, the same statement is true of the
algebras $A_{01}' \cap A_{i0}$.
Specifically there is some $i$ with
${\rm dim} Z(A_{01}' \cap A_{i0})>2L^2/\epsilon$.

For any $k>1$, we consider the inclusion
$A_{01}' \cap A_{i0} \subset A_{0k}' \cap A_{ik}$.  There must
be a minimal central projection $q_k$ of $A_{0k}' \cap A_{ik}$
such that ${\rm dim} Z(q_k(A_{01}' \cap A_{i0})) > 2/\epsilon$.
$q_k(A_{0k}' \cap A_{ik})$ is a factor, so we
may choose a projection $p_k<q_k$
in each 
$A_{0k}' \cap A_{ik}$ such that
$||E_{A_{01}' \cap A_{i0}}(p_k)||<\epsilon/2$.

For any $y \in A_{00}' \cap A_{i0}$, from Lemma 3 and Ocneanu compactness
$$
\lim_{k \rightarrow \infty}
\tr(E_{A_{0k}' \cap A_{ik}}(y) p_k) - \tr(E_{A_{01}' \cap A_{i0}}(y)p_k)=0$$
$$=\lim_{k \rightarrow \infty}
\tr(y p_k) - \tr(y E_{A_{01}' \cap A_{i0}}(p_k))$$
Since
$||E_{A_{01}' \cap A_{i0}}(p_k)||_2 < ||E_{A_{01}' \cap A_{i0}}(p_k)||<
\epsilon/2$, this means that $|\tr(y p_k)| < \epsilon||y||_2$
for $k$ sufficiently large.  $A_{00}' \cap A_{i0}$ is finite dimensional;
let its dimension be $d$.
By picking an orthonormal basis for this vector space,
we may find $k_0$ such that
$|\tr(x p_k)|<d\epsilon ||x||_2$ for all $x \in A_{00}' \cap A_{i0}$, $k >k_0$.
This means that $||E_{A_{00}' \cap A_{i0}}(p_k)||_2<d\epsilon$ for such $k$.
Taking $\tau$ to be the smallest trace of a nonzero projection in
$A_{00}' \cap A_{i0}$, we may bound the operator norm of
$E_{A_{00}' \cap A_{i0}}(p_k)$ as well:
$||E_{A_{00}' \cap A_{i0}}(p_k)|| < d \tau^{-1} \epsilon$ for
$k > k_0$.  $d$ and $\tau$ do not depend on $k$.

From
Lemma 4, we then have
$[(M_{2k})_{\omega} : M_0' \cap M_1^{\omega}] > \epsilon \tau/d$
for all $k > k_0$.  $\epsilon>0$ is arbitrary, and $d$ and $\tau$
are fixed,
so this means that
$\lim_{k \rightarrow \infty}
[(M_{2k})_{\omega} : M_0' \cap M_1^{\omega}]= \infty$.
\end{proof}

The lemmas from section 4.1 allow us to show that
$[(M_1)_{\omega}:M_0' \cap M_1^{\omega}]
\ge
[(M_k)_{\omega}:M_0' \cap M_1^{\omega}]$.  Combining this
fact with
Lemma 10 above gives the main result of this paper.

\begin{thm}
Let $M_0 \subset M_1$ be an commuting-square subfactor of infinite
depth.  Then the induced central sequence subalgebra
$M_0^{\omega} \cap M_1' \subset (M_1)_{\omega}$ has infinite index.
\end{thm}
\begin{proof}
Let the subfactor $M_0 \subset M_1$ have Jones tower
$M_0 \subset M_1 \subset M_2 \subset ...$

We choose $k \in \mathbb{N}$, and
consider the quadrilateral of Von Neumann algebras
$$\commsq{M_0^{\omega} \cap M_1'}{(M_1)_{\omega}}{(M_k)_{\omega}}
{M_1' \cap M_k^{\omega}}$$
Since $M_1 \subset M_k$ is a hyperfinite subfactor,
from Lemma 8 we know that
$[M_1' \cap M_k^{\omega}:(M_1)_{\omega}] = [M_k:M]$.
$M_1 \subset M_k$ is itself a commuting-square subfactor from
section 3.2, and so is extremal by Theorem 1.  This means
that Lemma 9 applies, and $[M_1' \cap M_k^{\omega}:(M_k)_{\omega}] \ge
[M_k:M_1]$.

From Lemma 5, all the algebras in the above quadrilateral are factors,
except possibly $M_0^{\omega} \cap M_1'$ in the lower left.  We have
$[M_k^{\omega} \cap M_1':(M_1)_{\omega}] < \infty$
and $[M_k^{\omega} \cap M_1':M_1^{\omega}] \le
[M_k^{\omega} \cap M_1':(M_k)_{\omega}]$, so the conditions of
Lemma 7 are satisfied and we have
$[(M_1)_{\omega}:M_0^{\omega} \cap M_1'] \ge
[(M_k)_{\omega}:M_0^{\omega} \cap M_1']$ for all $k$.

The original commuting-square subfactor $M_0 \subset M_1$ is of infinite
depth, so by Lemma 10 the sequence of Pimsner-Popa indices
$([(M_k)_{\omega}:M_0^{\omega} \cap M_1'])$ goes to infinity
with $k$.  Therefore the index of the central sequence
subalgebra $M_0^{\omega} \cap M_1' \subset M_1^{\omega}$ must be
infinite.
\end{proof}

\bibliographystyle{ams}
\providecommand{\bysame}{\leavevmode\hbox to3em{\hrulefill}\thinspace}
\providecommand{\MR}{\relax\ifhmode\unskip\space\fi MR }
\providecommand{\MRhref}[2]{%
  \href{http://www.ams.org/mathscinet-getitem?mr=#1}{#2}
}
\providecommand{\href}[2]{#2}

\end{document}